\documentclass[12 pt]{amsart}
\usepackage{amsfonts}
\usepackage{amsmath}
\usepackage{amsthm}
\newtheorem{theorem}{Theorem}[section]
\newtheorem{corollary}[theorem]{Corollary}
\newtheorem{conjecture}[theorem]{Conjecture}
\newtheorem{lemma}[theorem]{Lemma}
\newtheorem{proposition}[theorem]{Proposition}
\theoremstyle{remark}
\newtheorem{remark}[theorem]{Remark}

\theoremstyle{definition}
\newtheorem{definition}[theorem]{Definition}

\numberwithin{equation}{section}

\newcommand{\Z}{\mathbb{Z}}
\newcommand{\C}{\mathbb{C}}
\newcommand{\Q}{\mathbb{Q}}
\newcommand{\R}{\mathbb{R}}
\newcommand{\F}{\mathbb{F}}

\newcommand{\h}{\mathfrak{h}}
\newcommand{\al}{\alpha}
\newcommand{\be}{\beta}

\newcommand{\om}{\omega}

\newcommand{\ep}{\epsilon}

\newcommand{\la}{\langle}

\newcommand{\ra}{\rangle}
\newcommand{\Vir}{\mathrm{Vir}}

\newcommand{\Hom}{{\rm Hom}}

\newcommand{\BW}{{\rm BW}_{16}}

\newcommand{\CC}{{\mathcal C}}

\newcommand{\VA}{\rm VA}
\newcommand{\M}{\mathbb{M}}

\newcommand{\NO}{\,{\raise0.25em\hbox{$\mathop{\hphantom{\cdot}}%
\limits^{_{\circ}}_{^{\circ}}$}}\,}
\newcommand{\aut}{\mathrm{Aut}\,}

\begin{document}

\title[Ising vectors in $V_{\Lambda}^+$]{Ising vectors in the vertex operator algebra
$V_{\Lambda}^+$ associated with the Leech lattice $\Lambda$}

\author[C.H. Lam]{Ching Hung Lam }%
\address[C.H. Lam]{Department of Mathematics and National Center for Theoretical Sciences,
National Cheng Kung University,Tainan, Taiwan 701}%
\email{chlam@mail.ncku.edu.tw}%

\author[H. Shimakura]{Hiroki Shimakura}%
\address[H. Shimakura]{Graduate School of Science, Chiba
University, Chiba
 263-8522, Japan}%
\email{shima@math.s.chiba-u.ac.jp}%

\thanks{C.\,H. Lam is partially supported by NSC grant 95-2115-M-006-013-MY2 of Taiwan}%
\thanks{H.\,Shimakura is partially supported by the Japan Society for the Promotion of
Science Research Fellowships for Young Scientists} %

\begin{abstract}
In this article, we study the Ising vectors in the vertex operator
algebra $V_\Lambda^+$ associated with the Leech lattice $\Lambda$.
The main result is a characterization of  the Ising vectors in
$V_\Lambda^+$. We show that for any Ising vector $e$ in $V_\Lambda^+$, there is a sublattice $E\cong
\sqrt{2}E_8$ of $\Lambda$ such that $e\in V_E^+$. Some properties about their
corresponding $\tau$-involutions in the moonshine vertex operator
algebra $V^\natural$ are also discussed. We show that there is no
Ising vector of $\sigma$-type in $V^\natural$.
Moreover, we compute the centralizer $C_{\aut V^\natural}(z,
\tau_e)$ for any Ising vector $e\in V_\Lambda^+$, where $z$ is a
$2B$ element in $\aut V^\natural$ which fixes $V_\Lambda^+$.
Based on this result, we also obtain an explanation
for the $1A$ case of an observation by Glauberman-Norton (2001), which
describes some mysterious relations between the centralizer of $z$
and some $2A$ elements commuting $z$ in the Monster and the Weyl
groups of certain sublattices of the root lattice of type $E_8$ .
\end{abstract}
\maketitle

\section{Introduction}

The study of vertex operator algebras (VOAs) as modules of the rational
Virasoro VOA $L\left( \frac{1}{2},0\right) $ was first initiated by
Dong, Mason and Zhu\thinspace \cite{dmz}.
Partially motivated by their work, Miyamoto\thinspace \cite{m1}
introduced the notion of an Ising vector (i.e., a Virasoro vector of
weight $2$ which generates a copy of the rational Virasoro VOA
$L\left( \frac{1}{2},0\right) $ inside a VOA). He also developed
a simple method to construct involutions of the automorphism group
from Ising vectors. Such an automorphism is often called a Miyamoto
involution or a $\tau$-involution and is very useful for studying
the automorphism groups of VOAs. Therefore, it is really natural to
characterize all Ising vectors in VOAs. For example, all Ising
vectors in VOAs associated with binary codes and the VOAs
$V_{\sqrt2R}^+$ associated with root lattices $R$ were described in
\cite{lsy}. It was also conjectured that if $L$ is an even lattice
without roots, then any Ising vector in $V_L^+$  is contained in a
subVOA $V_U^+$ for a sublattice $U$ of $L$ isomorphic to $\sqrt2A_1$
or $\sqrt2E_8$.

The moonshine VOA $V^\natural$ constructed by
Frenkel-Lepowsky-Meurman\,\cite{flm} is one of the most important
examples of VOAs, which can be written as
\[
V^\natural = V_\Lambda^+\oplus V_\Lambda^{T,+},
\]
where $V_\Lambda^+$ is the fixed point subVOA of the Leech lattice
VOA $V_\Lambda$ by an automorphism $\theta$ induced from the
isometry $\be \mapsto -\be, \be \in \Lambda$, $V^T_\Lambda$ is the
unique irreducible $\theta$-twisted module of
$V_\Lambda$ and $V_\Lambda^{T,+}$ is the fixed point subspace of
$V^T_\Lambda$ by an automorphism induced from $\theta$. By the
construction, there is a natural involution $z\in \aut V^\natural$
such that $z|_{V_\Lambda^+}=1$ and $z|_{V_\Lambda^{T,+}}=-1$ and
this automorphism $z$ belongs to the $2B$ conjugacy class of the
Monster group.

On the moonshine VOA $V^\natural$, Miyamoto showed that $\tau$-involutions correspond to the $2A$ involutions of the Monster
group and that there is a one to one correspondence between
$2A$-involutions of the Monster group and Ising vectors in $V^\natural$ by
using the results in \cite{co}. This gives an approach to explain some mysterious phenomena associated with $2A$-involutions of the Monster group
by using the theory of VOAs. For example, the McKay observation on
the affine $E_8$-diagram has been studied in \cite{lyy} using
Miyamoto involutions. In this perspective, the explicit description
of Ising vectors in $V^\natural$ may provide some new insights on
these phenomena.

\medskip

In this article, we shall study the Ising vectors in the VOA
$V_\Lambda^+$. The main result is a characterization of all Ising
vectors in $V_\Lambda^+$. In particular, the conjecture in
\cite{lsy} on a characterization of Ising vectors in $V_L^+$ holds
when $L$ is isomorphic to the Leech lattice $\Lambda$. In addition,
some properties about the corresponding $\tau$-involutions of these
Ising vectors in the moonshine VOA $V^\natural$ will be discussed.
We also show that there is no Ising vector of $\sigma$-type in
$V^\natural$. In other words, the involution $\tau_e$ is non-trivial
for any Ising vector $e\in V^\natural$. Moreover, we compute the
centralizer $C_{\aut V^\natural}(z, \tau_e)$ of $z$ and the
$\tau$-involution $\tau_e$ in the automorphism group of $V^\natural$
for an Ising vector $e\in V_{\Lambda}^+$. It comes out that if $e$
is obtained from some sublattice $E\cong \sqrt{2}E_8$ in the Leech
lattice, then the centralizer $C_{\aut V^\natural}(z, \tau_e)$ also
stabilizes the subVOA $V_E^+$ and $C_{\aut V^\natural}(\tau_e,z)$
acts on $V_E^+$ as the orthogonal simple group $\Omega^+(8,2)$.

It is well-known that the products of any two $2A$-involutions of
the Monster simple group $\M$ fall into one of the following nine
conjugacy classes \cite{co,atlas}:
\[
1A,\ 2A,\ 3A,\ 4A,\ 5A,\ 6A,\ 4B,\ 2B,\text{ or }  3C.
\]

John McKay observed there is an interesting correspondence between
these nine conjugacy classes of  $\M$ and the nine nodes of the
extended Dynkin diagram $\hat{E}_8$ as follows:

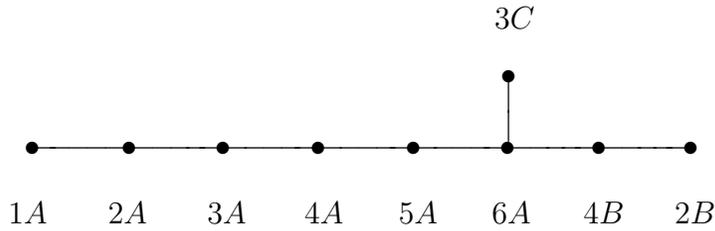
\begin{equation*}\label{E8}
\begin{array}{l}
  \hspace{184pt} 3C\\
  \hspace{186.2pt}\bullet \vspace{-13.1pt}\\
 \hspace{187.6pt}| \vspace{-10pt}\\
 \hspace{187.6pt}| \vspace{-13pt}\\
  \hspace{6pt} \bullet\hspace{-5pt}-\hspace{-7pt}-\hspace{-5pt}-\hspace{-5pt}-
  \hspace{-5pt}-\hspace{-5pt}\bullet\hspace{-5pt}-\hspace{-5pt}-
  \hspace{-5pt}-\hspace{-6pt}-\hspace{-7pt}-\hspace{-5pt}\bullet
  \hspace{-5pt}-\hspace{-5.5pt}-\hspace{-5pt}-\hspace{-5pt}-
  \hspace{-7pt}-\hspace{-5pt}\bullet\hspace{-5pt}-\hspace{-5.5pt}-
  \hspace{-5pt}-\hspace{-5pt}-\hspace{-7pt}-\hspace{-5pt}\bullet
  \hspace{-5pt}-\hspace{-6pt}-\hspace{-5pt}-\hspace{-5pt}-
  \hspace{-7pt}-\hspace{-5pt}\bullet\hspace{-5pt}-\hspace{-5pt}-
  \hspace{-6pt}-\hspace{-6pt}-\hspace{-7pt}-\hspace{-5pt}\bullet
  \hspace{-5pt}-\hspace{-5pt}-
  \hspace{-6pt}-\hspace{-6pt}-\hspace{-7pt}-\hspace{-5pt}\bullet
  \vspace{-6.2pt}\\
  \vspace{-10pt} \\
  1A\hspace{23pt} 2A\hspace{23  pt} 3A\hspace{22pt} 4A\hspace{21pt} 5A\hspace{20.5pt}
  6A\hspace{20pt} 4B\hspace{19pt} 2B\\
\end{array}
\end{equation*}
In Glauberman-Norton\,\cite{gn}, several other patterns and
mysterious properties related to the Dynkin diagram were discussed.
Among other things, they noted that the centralizer of a certain
subgroup generated by two $2A$-involutions and one $2B$-involution
in the Monster simple group appears to have a quotient which is
isomorphic to the ``half" of the Weyl group of the sub-diagram of
$\hat{E}_8$ with the relevant node removed.

In our case, the $2A$-involutions are $\tau_{e_E}$ and $\tau_{e_E}$
and the $2B$-involution is $z$. Then, $\tau_e\tau_e=1$ is of the
class $1A$. The subdiagram with the $1A$ node removed is of the type
$E_8$ and the Weyl group is of the shape $2.\Omega^+(8,2).2$. By our
method, the centralizer $C_{\aut V^\natural}(z, \tau_e)$ stabilizes
the subVOA $V_E^+$ and acts on $V_E^+$ as $\Omega^+(8,2)$, which is
the quotient of the commutator subgroup of the Weyl group of $E_8$
by its center. It explains the $1A$ case of the observation by
Glauberman and Norton. We believe that this is a first step toward a
complete understanding of the observation of Glauberman and Norton.
\medskip

\paragraph{\bf Acknowledgements.} Part of the work was done when the second author was visiting the National Center for Theoretical Sciences, Taiwan on January 2007.
He thanks the staff of the center for their help.

\medskip

\paragraph{\bf Notation.} Let $\Omega$ denote the set $\{1,2,\dots,24\}$. We
view the power set $\mathcal{P}(\Omega)$ of $\Omega$ as a
$24$-dimensional vector space over $\F_2$ naturally. For a lattice
$L$, let $L(m)$ denote the set of vectors in $L$ of squared norm
$2m$, namely $L(m)=\{v\in L\mid \langle v,v\rangle=2m\}$.

\section{Virasoro VOA and Ising vectors}

In this section, we shall recall some basic facts about Virasoro VOAs
and Ising vectors. We shall also review the construction of Miyamoto automorphisms.

Let $\Vir=\left(\bigoplus_{n\in \mathbb{Z}}\mathbb{C}L_n\right)\oplus
\mathbb{C} \mathbf{c}$ be the Virasoro algebra. Then $L_n$'s satisfy
the famous commutator relation
\[
[L_m,L_n]=(m-n)L_{m+n}+\frac{1}{12}(m^3-m)\delta_{m+n,0}\mathbf{c} .
\]
Let  $L\left( c,h\right) $ be the unique irreducible highest weight
module of $\Vir$ with central charge $c$ and highest weight $h$.
That means $L\left( c,h\right)$ is generated by a vector $v$ such
that $L_nv=0$ for $n>0$, $L_0v=hv$ and $\mathbf{c}v=cv$. It was
shown by Frenkel-Zhu\,\cite{fz} that $L(c,0)$ is a simple vertex
operator algebra. If $c=c_m=1-\frac{6}{(m+2)(m+3)}, m=1,2,3,\dots$,
then $L(c_m,0)$ has a unitary form and it is a rational VOA, that
means $L(c_m,0)$ has only finitely many irreducible modules and all
$L(c_m,0)$-modules are completely reducible (cf. \cite{dmz,w}).

When $c=c_1=\frac{1}{2}$, the VOA $L(\frac{1}{2},0)$ has exactly three inequivalent
irreducible modules, namely, $L\left( \frac{1}{2},0\right) ,L\left(
\frac{1}{2},\frac{1}{2}\right) $ and $L\left( \frac{1}{2},\frac{1}{16}\right) $, and its fusion rules
are given as
\begin{align*}
L\left( \frac{1}{2},\frac{1}{2}\right) \times L\left( \frac{1}{2},\frac{1}{2}%
\right) & =L\left( \frac{1}{2},0\right) , \\
L\left( \frac{1}{2},\frac{1}{2}\right) \times L\left( \frac{1}{2},\frac{1}{16%
}\right) & =L\left( \frac{1}{2},\frac{1}{16}\right) , \\
L\left( \frac{1}{2},\frac{1}{16}\right) \times L\left( \frac{1}{2},\frac{1}{%
16}\right) & =L\left( \frac{1}{2},0\right) +L\left( \frac{1}{2},\frac{1}{2}%
\right) \text{ }
\end{align*}
and $L\left( \frac{1}{2},0\right) $ acts as an identity.

\begin{definition}
A weight $2$ element $e\in V_2$ is called an \textsl{Ising vector} if the
vertex subalgebra $\VA(e)$ generated by $e$ is isomorphic to the simple
Virasoro VOA $L(\frac{1}{2},0)$ and $e$ is the Virasoro element of $\VA(e)$.
\end{definition}

\medskip

Next we shall review the definition of Miyamoto involutions. Let $e$
be an Ising vector. Since $\VA(e)\cong L(\frac{1}{2},0)$ is rational, we have the
isotypical decomposition
\[
V=V_e(0)\oplus V_e(\frac{1}2)\oplus V_e(\frac{1}{16}),
\]
where $V_e(h)$ denotes the sum of all irreducible
$\VA(e)$-submodules of $V$ isomorphic to $L(\frac{1}{2},h)$, $h=0,\frac{1}{2},
\frac{1}{16}$. An Ising vector $e$ is said to be \textsl{of $\sigma$-type} if $V_e(\frac{1}{16})=0$.

The following result was first proved by Miyamoto\,\cite{m1}:

\begin{theorem}
Let $V$ be a VOA and $e\in V$ an Ising vector. Then the linear map
$\tau_e:V\to V $ defined by
\[
\tau_e=
\begin{cases}
1& \text{ on } V_e(0)\oplus V_e(\frac{1}2)\\
-1& \text{ on } V_e(\frac{1}{16})
\end{cases}
\]
is an automorphism of $V$.

If $\tau_e={\rm id}$, namely $e$ is of $\sigma$-type, then the linear map $\sigma_e:V\to V$ defined by
\[
\sigma_e=
\begin{cases}
1& \text{ on } V_e(0)\\
-1& \text{ on } V_e(\frac{1}{2})
\end{cases}
\]
is an automorphism of $V$.
\end{theorem}

\begin{remark}
Note that the automorphism $\tau_e$ can also be defined by
\[
\tau_e= \exp( 16\pi i e_1),
\]
where $\exp(x)=\sum_{n=0}^\infty \frac{x^n}{n!}$ and
$Y(e,z)=\sum_{n\in \Z} e_n z^{-n-1}, e_n\in \mathrm{End}\, V$.
\end{remark}

\section{Lattice VOA and its $\Z_2$-orbifold} \label{sec:3}
Our notation for lattice VOAs here is standard (cf. \cite{flm}). Let
$L$ be a positive definite even lattice with inner product
$\la\,\cdot\, , \,\cdot\,\ra$. Then the VOA $V_L$ associated with
$L$ is defined to be $M(1) \otimes \C\{L\}$. More precisely, let
$\mathfrak{h} = \C \otimes_{\Z} L$ be an abelian Lie algebra and
$\hat{\mathfrak{h}} = \mathfrak{h} \otimes \C [t,t^{-1}] \oplus \C
K$ its affine Lie algebra. Let ${\hat{\frak{h}}}^-=\frak{h}\otimes
t^{-1}\C[t^{-1}]$ and let $S(\hat{\frak{h}}^-)$ be the symmetric
algebra of $\hat{\frak{h}}^-$. Then $M(1) =S(\hat{\frak{h}}^-)=
\C[\alpha(n)\mid\alpha \in \mathfrak{h}, n < 0]\cdot 1$ is the
unique irreducible $\hat{\mathfrak{h}}$-module such that $\alpha(n)
\cdot 1 = 0$ for $\alpha \in \mathfrak{h}$, $n \ge 0$ and $K=1$,
where $\alpha(n) = \alpha \otimes t^n$. Moreover,
$\C\{L\}=\mathrm{span}_\C\{e^\al\mid \al\in L\}$ is a twisted group
algebra of the additive group $L$ such that $e^\al e^\be=(-1)^{\la
\al,\be \ra} e^\be e^\al$. Note that $M(1)$ is a subVOA of $V_L$.

The twisted group algebra $\C\{L\}$ can be described by using
central extension. Let $\langle\kappa\rangle$ be a cyclic group of order $2$
and
\[
1\longrightarrow \langle\kappa\rangle  \longrightarrow
\hat{L}\bar{\longrightarrow}  L\longrightarrow 1
\]
a central extension of $L$ by $\langle\kappa\rangle$ with the commutator map
$c_0(\al, \be) = \la \al, \be \ra \mod 2$ for any $\al, \be \in
L$. Let $ L\to \hat{L}, \al\mapsto e^\al$ be a section.
Then the twisted group algebra
\[
\C\{L\}\cong \C[\hat{L}]/(\kappa+1) = \mathrm{span}_{\C}\{e^\al\mid
\al\in L\},
\]
where $\C[\hat{L}]$ is the usual group algebra of the group
$\hat{L}$.

Let $\theta: \hat{L}\to \hat{L}$ be an automorphism of $\hat{L}$ defined by
$\theta(a)=a^{-1}\kappa^{\la \bar{a},\bar{a}\ra/2}.$ Then it induces an automorphism of $V_L$ by
\[
\theta( \al_1(-n_1)\cdots \al_k(-n_k) e^\al)= (-1)^{k}
\al_1(-n_1)\cdots \al_k(-n_k)\theta (e^\al).
\]

Let $V_L^+=\{v\in V_L\mid \theta (v)=v\}$ be the fixed subspace of
$\theta$ in $V_L$. Then  $V_L^+$ is a subVOA and is often called a
\textsl{$\Z_2$-orbifold}.

For any sublattice $E$ of $L$, let
$\C\{E\}=\mathrm{span}_\C\{e^\al\in \C\{L\} \mid \al\in E\}$ be a
subalgebra of $\C\{L\}$ and let $\frak{h}_E=\C\otimes_\Z E$ be a
subspace of $\frak{h}= \C\otimes_\Z L$. Then the subspace
$S(\hat{\frak{h}}_E^-)\otimes \C\{E\}$ forms a subVOA of $V_L$ and
it is isomorphic to the lattice VOA $V_E$.

\medskip

Next, we shall review the constructions of some Ising vectors in a
lattice VOA $V_L$. The action of the corresponding $\tau$-involutions
will also be discussed.
Recall that two elements
$e,f$ in a VOA are said to be \textsl{mutually orthogonal} if
$Y(e,z)f=Y(f,z)e=0$.

\begin{theorem}[cf. \cite{dmz}]
Let $\al\in L(2)$. Then the elements
$\omega^+(\al)$ and $\omega^-(\al)$ defined by
\[
\omega^\pm(\al)=\frac{1}{16} \al(-1)^2\cdot \mathbf{1}\pm
\frac{1}4 (e^\al+\theta (e^{\al}))
\]
are two mutually orthogonal Ising vectors. 
\end{theorem}

The following proposition is easy:

\begin{proposition}\label{prp:3.2}
Let $\alpha\in L(2)$.
Then the automorphisms $\tau_{\omega^+(\al)}$ and $\tau_{\omega^-(\al)}$ of $V_L$ are equal to $\varphi_\al$ given by
\[
\varphi_\al (u\otimes e^\be)= (-1)^ {\la \al, \beta \ra} u\otimes e^\be
\quad \text{ for } u\in M(1)\text { and }\be \in L.
\]
\end{proposition}

If $L$ contains a sublattice isomorphic to $\sqrt{2}E_8$, then there
is another class of Ising vectors in $V_L$.

\begin{theorem}[cf. \cite{dlmn}]
Let $E\subset L$ be a sublattice isomorphic to $\sqrt{2}E_8$. Then
the Virasoro element of $V_E$ is given by
\[
\omega'=\frac{1}{240} \sum_{\al\in E(2)} \al(-1)^2\cdot \mathbf{1}
\]
and the element $e_E$ defined by
\[
e_E= \frac{1}{16}\omega' + \frac{1}{32} \sum_{\al\in E(2)} (e^\al+\theta(e^\al))
\]
is an Ising vector.
\end{theorem}

\begin{remark}
Note that the Ising vectors $\om^\pm(\al)$ and $e_E$ are clearly
contained in the orbifold VOA $V_L^+$
\end{remark}

The following theorem can be found in \cite{lyy}.

\begin{theorem}\label{thm:3.5}
Let $E\cong \sqrt{2}E_8$. Then $\tau_{e_E}$ is equal to
$\theta$ in $\aut V_E$. In particular, $\tau_{e_E}={\rm id}$ on $V_{E}^+$.
\end{theorem}

Now let $L$ be a positive definite even lattice and let
\[
L^*=\{ \al\in \Q\otimes_\Z L\mid \la \al, \be\ra\in \Z \text{ for all
} \be\in L\} \] be the dual lattice of $L$.  For $x\in L^*$, define
a $\Z$-linear map $\varphi_x: L \to \Z_2$ by
\[
\varphi_x(y)=\la x,y\ra \quad \mod 2.
\]
Clearly the map
\[
\begin{split}
\varphi: L^* & \longrightarrow \mathrm{Hom}_\Z(L, \Z_2)\\
         x &\longmapsto \ \varphi_x
\end{split}
\]
is a surjective group homomorphism and $\ker \varphi=2 L^*$. Hence, we have
$\mathrm{Hom}_\Z (L, \Z_2)\cong L^*/2L^*$.
For $\alpha\in L^*$, $\varphi_\alpha$ induces an automorphism of $V_L$ given
by
\[
\varphi_\alpha(u\otimes e^\beta) =(-1)^{\la \alpha,\beta\ra} u\otimes e^\beta \quad
\text{ for } u\in M(1)\text{ and } \beta \in L.
\]

Now suppose $L=E\cong \sqrt{2}E_8$. Then $E^*= \frac{1}2 E$ and
$E^*/2E^*= \frac{1}2 E/E\cong E_8/ 2E_8$. For $x\in
E^*=\frac{1}2 E$,
$$\varphi_x (e_E)= \frac{1}{16}\omega' + \frac{1}{32}
\sum_{\al\in E(2)} (-1)^{\langle
x,\al\rangle}(e^\al+\theta(e^\al))$$ is also an Ising vector in $V_E^+$ and we
have $256( = 2^8)$ Ising vectors of this form. Since $\varphi_x$
commutes $\theta$, $\tau_{\varphi_x (e_E)}=\theta$ on $V_E$ and
$\tau_{\varphi_x (e_E)}={\rm id}$ on $V_{E}^+$ also.

\begin{remark}\label{rmk:3.6}
It was shown in \cite{cs} that $E_8/2E_8$ contains  one class
represented by $0$, $120$ classes represented by a pair of roots
$\pm \alpha$, and $135$ classes represented by 16 vectors forming a
$4$-frame of $E_8$, i.e., a subset $\{\pm \al_1,\dots, \pm
\al_8\}\subset E_8$ such that $\la \al_i, \al_j\ra=4\delta_{i,j}$,
$i,j=1, \dots, 8$.
\end{remark}

\section{Leech Lattice and automorphism group of $V_\Lambda^+$}
In this section, we will review the automorphism group of $\Lambda$ and $V_\Lambda^+$ for the Leech lattice $\Lambda$.
\subsection{Leech lattice and automorphisms}
In this subsection, we will review some basic properties of the Leech
lattice $\Lambda$  and its automorphism group $Co_0$. Our notation
mainly follows that of Conway-Sloane\,\cite{cs}.

Let $\CC\subset P(\Omega)$ be the extended binary Golay code of length $24$.
A codeword $\mathcal{O}\in \CC$ is called
an \textsl{octad} if $|\mathcal {O}|=8$ and a \textsl{dodecad} if $|\mathcal{O}|=12$.
The automorphism group of the
Golay code $\CC$ is the Mathieu group $M_{24}$ and it acts
transitively on octads. Note also that $\CC$ is
generated by octads and there are exactly $759$ octads in $\CC$.




Let $\{e_i\mid i\in\Omega\}$ be an orthogonal basis of $\R^{24}$ of squared norm $2$.

\begin{theorem}[\cite{cs}]\label{gen}
The Leech lattice $\Lambda $ is a lattice of rank $24$ generated by
the vectors:
\begin{eqnarray*}
&&\frac{1}{2}e_{X}\,\,;\quad X\text{ is a generator of the Golay code }%
\CC, \\ &&\frac{1}{4}e_{\Omega }-e_{1}\,\,, \\ &&e_{i}\pm
e_{j}\,,\text{ }i,j \in \Omega
\end{eqnarray*}
where $e_{X}=\sum_{i\in X}e_{i}$.
\end{theorem}

\begin{definition}
A set of vectors $\{\pm \al_1, \dots,\pm \al_{24}\} \subset \Lambda$
is called an \textsl{$n$-frame} of $\Lambda$ if $\la \al_i, \al_j\ra= n
\delta_{i,j}$ for all $i,j\in \{1, \dots,24\}$.
\end{definition}

For example, $\{\pm 2e_1, \dots, \pm 2e_{24}\}$ is an $8$-frame of
$\Lambda$ and $\{ \pm(e_{2i-1}\pm e_{2i})\mid i=1,\dots, 12\}$ is a
$4$-frame of $\Lambda$.

\begin{lemma}{\rm (cf.\ \cite{cs})}\label{4frame}
Any vector in $\Lambda(2)$ is contained in a $4$-frame of $\Lambda$.
\end{lemma}
\begin{proof}
Since $\{ \pm(e_{2i-1}\pm e_{2i})\mid i=1,\dots, 12\}$ is a
$4$-frame of $\Lambda$ and the automorphism group $Co_0$ of
$\Lambda$ acts transitively on $\Lambda(2)$, we have the desired
result.
\end{proof}

Note that the Leech lattice $\Lambda$ has no vectors of squared norm
$2$ and is generated by the vectors of squared norm $4$.






Let $\mathcal{F}=\{\pm 2\al_1, \dots, \pm 2\al_{24}\}$ be an $8$-frame of
$\Lambda$. Then  $\{\al_1, \al_2, \dots,
\al_{24}\}$ is  an orthogonal basis of $\mathbb{R}^{24}$ such that
$\la \al_i,\al_j\ra = 2\delta_{i,j}$.

Denote
\[
\mathcal{C}_{\mathcal{F}}=\left\{S\subset \Omega=\{1,2, \dots,
24\}\,\left |\, \frac{1}{2} \sum_{i\in S}\, \al_i \in
\Lambda\right.\right\}.
\]
It is well-known (cf. \cite{cs,g2}) that $\mathcal{C}_{\mathcal{F}}$
is isomorphic to the binary Golay code $\mathcal{C}$. The Leech
lattice $\Lambda$ will then be generated by the vectors of the form
\[
2\al_i,\qquad  \al_i\pm \al_j, \qquad \frac{1}2 \al_S=\frac{1}2
\sum_{i\in S} \al_i,\qquad \text{ and } \qquad  \frac{1}4
\al_{\Omega}-\al_i,
\]
where $i,j \in \Omega$ and $S\in \mathcal{C}$.

For any permutation $\pi$ of $\Omega$, $\pi$ defines an isometry of
$\R^{24}$ by $\pi(\al_i)=\al_{\pi i}$. If $\pi(
\mathcal{C}_{\mathcal{F}})=\mathcal{C}_{\mathcal{F}}$, then $\pi$ also defines an automorphism
of $\Lambda$. Let $S$ be a subset of $ \Omega$. Then we can also
define an isometry $\varepsilon_S^{\mathcal{F}}: \R^{24} \to
\R^{24}$ by $\varepsilon_S^{\mathcal{F}}(\al_i)= -\al_i$ if $i\in S$
and $\varepsilon_S^{\mathcal{F}}(\al_i)= \al_i$ if $i\notin S$.






The involutions in the Conway group $Co_0$ can be characterized as
follows:

\begin{theorem}[\cite{cs,g2}]\label{thm:4.2}
There are exactly $4$ conjugacy classes of involutions in $Co_0$.
They correspond to the involutions $\varepsilon_S^{\mathcal{F}}$,
where $\mathcal{F}$ is an $8$-frame and $S\in \mathcal{C}$ is an
octad, the complement of an octad, a dodecad, or the set $\Omega$.
Moreover, the sublattice $\{v\in\Lambda\mid
\varepsilon_S^{\mathcal{F}}(v)=-v\}$ is isomorphic to $\sqrt2E_8$,
${\rm BW}_{16}$, $\sqrt2D_{12}$ and $\Lambda$ respectively, where
${\rm BW}_{16}$ is the Barnes-Wall lattice of rank $16$.
\end{theorem}

Note that the center of the Conway group $Co_0$ is the cyclic group
$\langle\pm 1\rangle$. The quotient group $Co_1=Co_0/ \langle\pm 1\rangle$ is a simple group
and it has $3$ conjugacy classes of involutions, namely $2A$, $2B$,
and  $2C$ (cf. \cite{atlas}).

\begin{proposition}[\cite{atlas,cs}]\label{prop:4.9}
Let $g$ be an involution of $Co_1$ and $\tilde{g}$ a lift of $g$ in
$Co_0$.
\begin{enumerate}
\item If $g$ is of class $2A$, then
$\tilde{g}=\varepsilon_\mathcal{O}^{\mathcal{F}}$ or
$\varepsilon_{\Omega+\mathcal{O}}^{\mathcal{F}}$ for some octad
$\mathcal{O}$ and an $8$-frame $\mathcal{F}$.

\item If $g$ is of class $2B$, then
$\tilde{g}$ is of order $4$.

\item If $g$ is of class $2C$, then
$\tilde{g}=\varepsilon_\mathcal{S}^{\mathcal{F}}$ for some dodecad
$\mathcal{S}$ and an $8$-frame $\mathcal{F}$.
\end{enumerate}
\end{proposition}

Now, we shall consider the embeddings of the lattice $\sqrt{2}E_8$
into the Leech lattice $\Lambda$. First we shall construct some
sublattices isomorphic to $\sqrt{2}E_8$ in $\Lambda$.

Let $\mathcal{F}=\{\pm 2\al_1, \dots, \pm 2\al_{24}\}$ be an $8$-frame of
$\Lambda$. Then  $\{\al_1, \al_2, \dots,
\al_{24}\}$ is  an orthogonal basis of $\mathbb{R}^{24}$ such that
$\la \al_i,\al_j\ra = 2\delta_{i,j}.$

\begin{lemma}
Let $\mathcal{O}\in \mathcal{C}_\mathcal{F}$ be an octad and denote
\[
E_{\mathcal{F}}(\mathcal{O})=\mathrm{span}_\Z\left\{ \al_i\pm
\al_j, \frac{1}2 \sum_{i\in \mathcal{O}} \ep_i \al_i\right \},
\]
where $i,j\in \mathcal{O}$ and  $\ep_i=\pm 1$ such that $\prod_{i\in
\mathcal{O}} \ep_i=1$. Then $E_\mathcal{F}(\mathcal{O})$ is
isomorphic to $\sqrt{2}E_8$.
\end{lemma}

\begin{theorem}[\cite{g2}] The Conway group $Co_0$ is transitive on the
set $ \{ E\subset \Lambda\mid E\cong \sqrt{2}E_8\}$ of sublattices of the Leech lattice isomorphic to $\sqrt2E_8$.
\end{theorem}

As a corollary, we have the following proposition.

\begin{proposition}\label{prop:4.15}
Let $E\subset \Lambda$ be a sublattice isomorphic to $\sqrt{2}E_8$.
Then there exists an $8$-frame $\mathcal{F}$ of $\Lambda$ and an octad
$\mathcal{O}\in \CC_\mathcal{F}$ such that $E_{\mathcal{F}}(\mathcal{O}) =E$.
\end{proposition}

\begin{proof}
Fix an $8$-frame $\mathcal{F}_0=\{\pm \be_1, \cdots, \pm \be_{24}\}$
and an octad $\mathcal{O}\in \mathcal{C}_{\mathcal{F}_0}$. For any $E\cong
\sqrt{2}{E_8}$, there is $g\in \aut\, \Lambda$ such that $g(E)=
E_{\mathcal{F}_0}(\mathcal{O})$. Then $\mathcal{F}=\{ \pm g^{-1}(\be_1), \dots, \pm g^{-1}(\be_{24})\}$ forms an $8$-frame of $\Lambda$ and $E$ is clearly equal to $E_{\mathcal{F}}(\mathcal{O})$, where we regard $\mathcal{O}$ as an element in $\mathcal{C}_{\mathcal{F}}$.
\end{proof}

\begin{lemma}\label{rmk:4.14}{\rm (cf.\ \cite{cs,g2})}
For each $E\cong \sqrt{2}E_8\subset \Lambda$, there exists
$135$ distinct $8$-frames $\mathcal{F}$ of $\Lambda$ such that
$E=E_\mathcal{F}(\mathcal{O})$.
\end{lemma}
\begin{proof} Let $\mathcal{F}_0$ be an $8$-frame of $E$.
Then $\mathcal{F}_0$ forms a coset of $E/2E$ and, $\mathcal{F}_0$ is contained in a coset
of $\Lambda/2\Lambda$. This implies that there exists a unique
$8$-frame $\mathcal{F}$ of $\Lambda$ such that
$\mathcal{F}_0\subset\mathcal{F}$ (cf.\ \cite{cs,g2}). It follows
from Remark \ref{rmk:3.6} that $E$ contains exactly $135$ distinct
$8$-frames, which proves this lemma.
\end{proof}

\subsection{Automorphism group of $V_\Lambda^+$}
The full automorphism group of $V_\Lambda^+$ has been determined in
\cite{sh}. We shall recall its main results.

Let $L$ be a positive definite even lattice and
\[
1 \longrightarrow \langle\kappa\rangle \longrightarrow \hat{L} \bar{
\longrightarrow} L \longrightarrow 1
\]
a central extension of $L$ by $\langle\kappa\rangle\cong\Z_2$ such that the commutator map $c_0(\al, \be)={\la \al, \be\ra} \mod 2$, $\al,\be\in L$.
The following theorem is well-known (cf. \cite{flm}):

\begin{theorem}\label{thm:4.10} For an even lattice $L$, the sequence
\[
1 \longrightarrow \mathrm{Hom}(L, \Z_2) { \longrightarrow}
\aut \hat{L} \longrightarrow \aut L\longrightarrow  1
\]
is exact.
In particular $\aut \hat{\Lambda}\cong 2^{24} \cdot Co_0.$
\end{theorem}

Recall that $\theta$ is the automorphism of $V_\Lambda$ defined by
\[
\theta( \al_1(-n_1)\cdots \al_k(-n_k) e^\al)= (-1)^{k}
\al_1(-n_1)\cdots \al_k(-n_k)\theta (e^\al),
\]
where $\theta(a)=a^{-1}\kappa^{\la \bar{a},\bar{a}\ra/2}$ on
$\hat{L}$.

\begin{lemma}([\cite{sh}])\label{lem:4.11}
Let $L$ be a positive definite even lattice without roots, i.e.,
$L(1)=\emptyset$. Then the centralizer
$C_{\aut{V_L}}(\theta)$ of $\theta$ in $\aut{V_L}$ is isomorphic to
$\aut{\hat{L}}$. If $L=\Lambda$ is the Leech lattice, we have
\[
C_{\aut{V_\Lambda}}(\theta)\cong \aut{\hat{\Lambda}}\cong
2^{24}\cdot Co_0.
\]
\end{lemma}

\begin{proposition}\label{prp:4.17}
Let $E= E_\mathcal{F}(\mathcal{O})\cong {\sqrt{2}E_8}$ be a
sublattice of $\Lambda$ and $ e$ an Ising  vector in $V_\Lambda$ of
the form $\varphi_x(e_E)$. Then $\tau_e\in
C_{\aut{V_\Lambda}}(\theta)$ and the image of $\tau_e$ under the
canonical epimorphism from $C_{\aut{V_\Lambda}}(\theta)(\cong
\aut \hat{\Lambda})$ to
$C_{\aut{V_\Lambda}}(\theta)/\Hom(\Lambda,\Z_2)\cong
\aut \Lambda =Co_0$ is
$\varepsilon_\mathcal{O}^{\mathcal{F}}$.
\end{proposition}

\begin{proof}
Since $e$ is fixed by $\theta$, it is clear that $\tau_e\in
C_{\aut{V_\Lambda}}(\theta)$.

Let $g$ be the image of $\tau_e$ under the canonical
epimorphism from $\aut \hat{\Lambda}$ to $\aut \Lambda$. Let
$B=\{\al\in \Lambda\mid \la \al, \be\ra=0 \text{ for all } \be \in
E\}$ be the orthogonal complement of $E$ in $\Lambda$. Then $B$ is
isomorphic to the Barnes-Wall lattice ${\rm BW}_{16}$ (cf. \cite{cs,g2}).
By the definition of $e$, $e_1v=0$ for all $v\in V_B$. Hence
${\tau_e}_{|_{V_B}}={\rm id}$ but $\tau_e$ acts as $\theta$ on $V_E$
by Theorem \ref{thm:3.5}. Thus $g$ acts as $1$ on $B$ and
$-1$ on $E$. As $E=E_\mathcal{F}(\mathcal{O})$, we have
$g=\varepsilon_\mathcal{O}^{\mathcal{F}}$ as desired.
\end{proof}

\begin{theorem}([\cite{sh}])
Let $V_\Lambda^+=\{ v\in V_\Lambda\mid \theta(v)=v\}$ be the fixed
point subVOA of $\theta$ in $V_\Lambda$. Then
$
\aut{V_\Lambda^+}\cong C_{\aut{V_\Lambda}}(\theta)/ \langle\theta\rangle \cong
2^{24}\cdot Co_1
$
and the sequence
\[
1 \longrightarrow  \mathrm{Hom}(\Lambda, \Z_2) \longrightarrow
\aut{V_\Lambda^+}\ \overset{\xi}{\longrightarrow} \ \aut \Lambda/\langle\pm
1\rangle\longrightarrow  1.
\]
is exact.
\end{theorem}

\section{Ising Vectors in the $\Z_2$-orbifold VOA $V_\Lambda^+$}
In this section, we shall describe the Ising vectors in
$V_\Lambda^+$ and discuss the corresponding symmetries of the
moonshine VOA $V^\natural$.

\subsection{Characterization of Ising vectors in $V_\Lambda^+$.}

In \cite{lsy}, all Ising vectors in the VOA $V_{\sqrt{2}R}^+$ has
been characterized for a root lattice $R$ as follows.
\begin{proposition}({\rm \cite{lsy}})\label{lsy2} Let $R$ be a root lattice.
Then any Ising vector in $V_{\sqrt{2}R}^+$ is either equal to $\om ^\pm(\al)$ for some $\al\in \sqrt2R$ of norm $4$ or
$\varphi_x(e_E)$ for some sublattice $E$ of $\sqrt2R$ isomorphic to $\sqrt{2}E_8$ and $x\in \frac{1}{2}E$.
\end{proposition}

Note that $\om^\pm(\al)\in V_{\Z\al}^+\cong V_{\sqrt2A_1}^+$ and that $\varphi_x(e_E)\in V_{E}^+\cong  V_{\sqrt2E_8}^+$.
Its generalization was given in \cite{lsy} as a conjecture.

\begin{conjecture} Let $L$ be an even lattice $L$ without roots. Then any Ising
vector of $V_L^+$ belongs to a subVOA $V_U^+$ for some sublattice
$U$ of $L$ isomorphic to $\sqrt2A_1$ or $\sqrt2E_8$.
\end{conjecture}

We shall show that the conjecture holds when $L$ is isomorphic to
the Leech lattice $\Lambda$. First let us recall a result from
\cite{lsy}.

\begin{lemma}\label{LLSY}{\rm \cite[Lemma 3.7]{lsy}} Let $V=\oplus_{i=0}^\infty V_n$ be a VOA with
$V_0=\C{\bf 1}$, $V_1=0$. Suppose that $V$ has two Ising vectors $e$, $f$ and that
$e$ is of $\sigma$-type, i.e., $\tau_e=\mathrm{id}$ on $V$. Then $e$
is fixed by $\tau_f $, namely $e\in V^{\tau_f}$.
\end{lemma}

Let us discuss Ising vectors of $V_L^+$ and the associated Miyamoto
involutions.

\begin{lemma}\label{L12} Let $L$ be an even lattice of rank $n$ without roots.
Let $U$ be the sublattice of $L$ generated by $L(2)$. Suppose that
$U\cong\sqrt2R$ for some root lattice $R$ of rank $m$. Then any
Ising vector of $V_L^+$ belongs to $V_U^+$.
\end{lemma}
\begin{proof}
Let $U^\perp=\{\al\in L\mid \la \al, \be\ra=0 \text{ for all } \be \in
U\}$ be the orthogonal complement of $U$ in $L$. Let
$\mathfrak{h}_0=\C\otimes_\Z U$ and $\mathfrak{h}_1=\C\otimes_\Z
U^\perp$. Then $\mathfrak{h}=\C\otimes_\Z L=\mathfrak{h}_0\perp
\mathfrak{h}_1$. Since  the weight $2$ subspace of $(V_L^+)_2$ is
equal to the weight $2$ subspace of $[S(\hat{\mathfrak{h}}^-)\otimes
\C\{U\}]^+\cong [ V_U\otimes S(\hat{\h}_1^-)]^+$, we know that any
Ising vector $e$ of $V_L^+$ must be in $[ V_U\otimes
S(\hat{\h}_1^-)]^+$.

Let $U^\prime\cong(\sqrt2A_1)^{m-n}$ and $K=U\perp U^\prime$. Then
both $L$ and $K$ are of rank $n$ and we can identify $\C\otimes_\Z
K$ with $\mathfrak{h}$ and $\C\otimes_\Z U'$ with $\mathfrak{h}_1$.
Then the lattice VOA $V_K$ is given by
\[
V_K = S(\hat{\mathfrak{h}}^-)\otimes \C\{K\}\cong
S(\hat{\h}_0^-)\otimes S(\hat{\h}_1^-)\otimes \C\{U\}\otimes
\C\{U'\},
\]
which contains $V_U \otimes S(\hat{\h}_1^-)$ as a subalgebra. Hence,
we can view the Griess algebra of $V_{L}^+$ as a Griess subalgebra
of $V_K^+$, also. Thus by Proposition \ref{lsy2}, $e\in V_U^+$ or
$e\in V_{U^\prime}^+$ as $K/\sqrt2$ is isomorphic to a root lattice.
It follows that $e\in V_{U}^+$ as there is no Ising vector in
$S(\hat{\h}_1^-)$.

\end{proof}

\begin{lemma}\label{L13} Let $L$ be an even lattice without roots and $e$
an Ising vector of $V_L^+$.
Then the automorphism $\tau_e$ belongs to $C_{\aut{V_L}}(\theta)/\langle\theta\rangle$.
In particular, $\tau_e\in \aut{\hat{L}}/\langle\theta\rangle$.
\end{lemma}

\begin{proof}
We view $\tau_e$ as an automorphism of $V_L$. Since $\theta$ fixes
$e$, we have $\theta\tau_e\theta=\tau_{\theta(e)}=\tau_e$, which
proves this lemma.
\end{proof}

Since $C_{\aut{V_\Lambda}}(\theta)/\langle\theta \rangle \cong
2^{24}\cdot Co_1$, it suffices to consider the two cases:
$\xi({\tau}_e)=1$ and $\xi({\tau}_e)\neq 1$, where
$\xi:C_{\aut{V_\Lambda}}(\theta)/\la \theta \ra\to Co_1$ denotes
the natural epimorphism. First, we study the case $\xi({\tau}_e)=1$
for general $L$.

\begin{proposition}
Let $L$ be an even lattice without roots and let $e$ be an Ising vector of
$V_{L}^+$. Suppose that $\xi({\tau_e})=1$.
Then $e$ is equal to $\om^{\pm}(\alpha)$ for
some $\alpha\in L(2)$ or $\varphi_x(e_E)$ for some $E\subset L$ isomorphic to $\sqrt2E_8$ and $x\in \frac{1}{2}E$.
\end{proposition}

\proof\ By Lemma \ref{L13}, $\tau_e\in
C_{\aut{V_L}}(\theta)/\langle\theta\rangle$. By the assumption,
$\tau_e=\varphi_v$ for some $v\in L^*/2L^*$. Then $e\in
(V_{L}^{+})^{\varphi_v}$ and by Proposition \ref{prp:3.2}
$(V_{L}^+)^{\varphi_v}=V_{L_v}^+$, where $L_v=\{u\in L\mid \langle
v,u\rangle\in2\Z\}$. Note that $e$ is of $\sigma$-type in
$V_{L_v}^+$. Set $A=\langle \tau_{\om^\pm(\alpha)}\mid \alpha\in
L_v(2)\rangle$. By Lemma \ref{LLSY}, $e\in (V_{L_v}^{+})^{A}$.
Moreover $(V_{L_v}^{+})^{A}=V_{\tilde{L}_v}^+$, where $\tilde{L}_v=\{u\in
L_v\mid \langle u,w\rangle\in2\Z,\ \forall w\in L_v(2)\}$. Let $U$ be the
sublattice generated by $\tilde{L}_v(2)$.  In this case,
$(1/\sqrt2)U$ is isomorphic to a root lattice. Thus, by Lemma
\ref{L12}, $e\in V_U^+$ and $e=\om^{\pm}(\alpha)$ for some
$\alpha\in \tilde{L}_v(2)$ or $e=\varphi_x(e_E)$ by Proposition
\ref{lsy2}. \qed

Next, we consider the case $\xi({\tau}_e)\neq 1$ for general $L$.

\begin{lemma}\label{L15}Let $L$ be an even lattice without roots.
Let $e$ be an Ising vector in $V_L^+$ such that $\xi(\tau_e)\neq1$ and let $g\in
\aut L$ be a lift of $\xi({\tau_e})$. Then
$e\in V_{L\cap 2N(e)^*}^+$, where $N(e)$ is the sublattice generated by $\{ v\in L(2)\mid g(v)=\pm v\}$.
\end{lemma}

\begin{proof} Set $U^i=\{v\in L\mid g(v)=(-1)^i v\}$ ($i=0,1$) and
$U=U^0\oplus U^1$. Since for $\lambda\in L$ both
$\lambda+g(\lambda)$ and $\lambda-g(\lambda)$ belong to $U$, we have
$2L\subset U$. Consider the coset $L/U\subset ((U^0)^\perp\oplus
(U^1)^\perp)/(U^0\oplus U^1)$ and the canonical projections
$\lambda\in L/U$ to $\lambda_i\in (U^i)^\perp/U^i$ for $i=0.1$. Then
we obtain the decomposition
$$(V_L^{+})^{\tau_e}=\bigoplus_{\lambda_0+\lambda_1\in
L/U}V_{\lambda_0+U_0}^{\tau_e\theta}\otimes
V_{\lambda_1+U_1}^{\tau_e},$$ where
$V_{\lambda_0+U_0}^{\tau_e\theta}$ is the $\tau_e\theta$-fixed point
subspace of $V_{\lambda_0+U_0}$ and $V_{\lambda_1+U_1}^{\tau_e}$ is
the $\tau_e$-fixed point subspace of $V_{\lambda_1+U_1}$. Set
$A=\langle\tau_{\omega^\pm(\alpha)}\mid \alpha\in U(2)\rangle$.
Then $N(e)$ is a sublattice of $L$ generated by $U(2)$ and $A=\langle\varphi_\alpha\mid \alpha\in N(e)\rangle$.
Since $e$ is of $\sigma$-type in $(V_L^{+})^{\tau_e}$, $e\in
(V_L^{+})^{\langle\tau_e,A\rangle}=(V_{L\cap 2N(e)^*}^{+})^{\tau_e}\subset V_{L\cap
2N(e)^*}^+$ by Lemma \ref{LLSY}.\end{proof}

We apply the lemma above to characterize Ising vectors of
$V_{\Lambda}^+$. Let $e$ be an Ising vector in $V_{\Lambda}^+$ such
that $\xi({\tau}_e)\neq 1$. Then by Proposition \ref{prop:4.9}, the
conjugacy class of $\xi({\tau}_e)\in Co_1$ is $2A$ or $2C$ since
$\tau_e$ is an involution in $\aut{V_\Lambda}$. Hence, we should
consider these two cases.

\begin{proposition} Let $e$ be an Ising vector in $V_\Lambda^+$. Suppose that
$\xi({\tau_e})$ belongs to the conjugacy class $2A$ of $Co_1$. Then
$N(e)=E\oplus E^\prime$ and $e=\varphi_x(e_E)$, where
$E\cong\sqrt2E_8$, $E^\prime\cong\BW$ and $x\in \frac{1}{2}E$.
\end{proposition}

\begin{proof} By Proposition \ref{prop:4.9},
$\xi({\tau}_e)=\varepsilon_\mathcal{O}^{\mathcal{F}}$ or
$\varepsilon_{\Omega+\mathcal{O}}^{\mathcal{F}}$ for some octad
$\mathcal{O}$. Thus, by Theorem \ref{thm:4.2}, $N(e)=E\oplus
E^\prime$, where $E\cong\sqrt2E_8$ and $E^\prime \cong\BW$. Then
$\Lambda\cap 2N(e)^*\cong \sqrt2E_8\oplus2\BW^*$. By Lemma
\ref{L15}, $e\in V_{\sqrt2E_8\oplus2\BW^*}^+$. Since
$2\BW^*\cong\sqrt2\BW$, we have $2\BW^*(2)=\emptyset$. Then by Lemma
\ref{L12}, $e\in V_{E}^+$. Hence by Proposition \ref{lsy2},
$e=\omega^\pm(\alpha)$ or $\varphi_x(e_E)$. By Proposition
\ref{prp:3.2}, $e$ must be $\varphi_x(e_E)$ for some $x\in \frac{1}{2}E$ as
$\xi({\tau}_e)\neq 1$.
\end{proof}

\begin{proposition} There is no Ising vector $e$ in $V_\Lambda^+$ such that
$\xi({\tau_e})$ belongs to the conjugacy class $2C$ of $Co_1$.
\end{proposition}

\begin{proof} Let $e$ be an Ising vector in $V_\Lambda^+$ such that
$\xi({\tau_e})$ belongs to the conjugacy class $2C$ of $Co_1$.  By
Theorem \ref{thm:4.2} and Proposition \ref{prop:4.9}, $N(e)\cong
\sqrt2D_{12}\oplus\sqrt2D_{12}$. Then $\Lambda\cap
2N(e)^*\cong\sqrt2D_{24}^+$. By Lemma \ref{L15}, $e\in
V_{\sqrt2D_{24}^+}^+$. Since the set of norm $4$ vectors in
$\sqrt2D_{24}^+$ generates $\sqrt2D_{24}$, $e$ belongs to
$V_{\sqrt2D_{24}}^+$. Note that there is no sublattice of
$\sqrt2D_{24}$ isomorphic to $\sqrt2E_8$. Hence, by Proposition
\ref{lsy2}, $e=\om^\pm(\al)$ for some $\al\in \sqrt{2}D_{24}(2)$. In
this case, $\xi({\tau}_e)=1$ by Proposition \ref{prp:3.2}, which is
a contradiction.
\end{proof}

Therefore, the conjecture holds when $L$ is isomorphic to the Leech
lattice. In particular, there are exactly $2$ types of Ising vectors
in $V_\Lambda^+$, namely, $\omega^\pm(\alpha)$ for some $\alpha\in
\Lambda(2)$ or $\varphi_x( e_E)$ for some sublattice $E\cong
\sqrt{2}E_8$ and $\frac{1}2x\in E$

Since $Co_0$ is transitive on $\Lambda(2)$, for any $\al\in
\Lambda(2)$, there is a sublattice $E\cong \sqrt{2}E_8$ such that
$\al\in E$. Therefore we also have
\begin{theorem}\label{MainTheorem}
Let $\Lambda$ be the Leech lattice. Then for any Ising vector $e$ in
$V_\Lambda^+$, there is a sublattice $E\cong \sqrt{2}E_8$ such that
$e\in V_E^+$.
\end{theorem}

As a corollary, we can count all Ising vectors in $V_\Lambda^+$.
\begin{corollary}
There are exactly $11935319760$ Ising vectors in the VOA
$V_\Lambda^+$.
\end{corollary}

\begin{proof}
Since $\Lambda$ has exactly $196560$ vectors of squared norm $4$,
there are $196560$ Ising vectors of the form $\om^\pm(\al), \al\in
\Lambda(2)$. By Proposition \ref{prop:4.15}, any sublattice isomorphic to $\sqrt2E_8$ is determined by an $8$-frame of $\Lambda$ and an octad.
By Lemma \ref{rmk:4.14}, there are exactly
\[
|\{\mathcal{O}\in\CC|\ |\mathcal{O}|=8\}|\times
\frac{|\Lambda(4)|}{48}\times\frac{1}{135}=759\times
\frac{398034000}{48}\times \frac{1}{135}= 46621575
\]
sublattices isomorphic to $\sqrt{2}E_8$. Thus, there are $11935123200
(\ =46621575\times 256\,)$ Ising vectors of the form $\varphi_x
(e_E)$. Therefore, there are totally $11935319760$ Ising vectors in
$V_{\Lambda}^+$.
\end{proof}

\begin{corollary}\label{cor}
None of the Ising vectors in $V_\Lambda^+$ is of $\sigma$-type,
namely $V_{e}(\frac{1}{16})\neq0$ for any Ising vector $e\in V_\Lambda^+$.
\end{corollary}
\begin{proof} Let $e$ be an Ising vector in $V_\Lambda^+$.
Then by Theorem \ref{MainTheorem}, $e=\omega^\pm(\alpha)$ or $e=\varphi_x( e_E)$, and $\tau_e$ is not of $\sigma$-type by Proposition \ref{prp:3.2} and \ref{prp:4.17} respectively.
\end{proof}

Let us consider the correspondence between Ising vectors
in $V_{\Lambda}^+$ and Miyamoto automorphisms.

\begin{proposition}\label{prp:5.11} Let $e,f$ be Ising vectors in $V_\Lambda^+$ such that $\tau_e=\tau_f\in\aut{V_{\Lambda}^+}$.
\begin{enumerate}
\item If $e=\omega^\pm(\alpha)$ then $f$ is equal to $\omega^+(\alpha)$ or $\omega^-(\alpha)$.
\item If $e=\varphi_x(e_E)$ for some sublattice $E\cong\sqrt2E_8$ of $\Lambda$ and $x\in \frac{1}{2}E$ then $f=e$.
\end{enumerate}
\end{proposition}
\begin{proof}
First we suppose that $e=\omega^\pm(\alpha)$.
Then by Proposition \ref{prp:3.2} and \ref{prp:4.17}, $f=\omega^\pm(\beta)$ for some $\beta\in \Lambda(2)$, and $\langle\beta, v\rangle=\langle\alpha,v\rangle\mod 2$ for all $v\in \Lambda$.
Since $\Lambda$ is unimodular, $\beta\in \alpha+2\Lambda$, and (1) follows.

Next we suppose that $e=\varphi_x(e_E)$. By Proposition
\ref{prp:3.2} and \ref{prp:4.17}, $f$ must be $\varphi_y(e_E)$ for
some $y\in \frac{1}{2}E$. Then $$\tau_e\tau_f=\varphi_{2(x+y)}={\rm id}$$ since $\tau_{\varphi_v(e_E)}=\varphi_{2v}\tau_{e_E}$ for $v\in \frac{1}{2}E$. Recall that if $u\in E$ satisfies $\langle u,\Lambda\rangle\in2\Z$ then $u\in 2E$. Hence
$2(x+y)\in 2E$, namely $y\in x+ E$, and we obtain $e=f$.
\end{proof}

At end of this subsection, we give some applications to Ising vectors
in the moonshine VOA $V^\natural=V_{\Lambda}^+\oplus
V_{\Lambda}^{T,+}$. Let $z$ be the linear map of $V^\natural$ acting
as $1$ and $-1$ on $V_{\Lambda}^+$ and $V_{\Lambda}^{T,+}$
respectively. Then $z$ is an automorphism of $V^\natural$.

\begin{lemma}\label{lem:5.12} Let $\alpha\in \Lambda(2)$.
Then $\tau_{\om^+(\al)} \tau_{\om^-(\al)}=z$ on $V^\natural$.
\end{lemma}

\begin{proof}
Set $g=\tau_{\om^+(\al)} \tau_{\om^-(\al)}$.
By Proposition \ref{prp:3.2}, $g=1$ on $V_\Lambda^+$.
Since the $V_\Lambda^+$-module $V_\Lambda^{T,+}$ is irreducible, $g$ must be $1$ or $-1$ on $V_{\Lambda}^{T,+}$.

Suppose that $g=1$ on
$V_{\Lambda}^{T,+}$. Let $F=\{\al_1=\al, \al_2, \dots, \al_{24}\}$
be a $4$-frame of $\Lambda$ containing $\alpha$ (cf.\ Lemma \ref{4frame}) and
$\mathcal{T}=\{\om^\pm(\al_1), \dots, \om^\pm(\al_{24})\}$ the
Virasoro frame of $V^\natural$ associated with $F$. Let $(C,D)$ be
the structure code of $V^\natural$ with respect to $\mathcal{T}$,
namely $C$ and $D$ are the binary codes of length $48$ such that
$V^0$ is a code VOA associated with $C$ and $V^\natural=\oplus_{d\in
D} V^d$, where $V^d$ is an irreducible $V^0$-module whose
$\frac{1}{16}$-word is $d\in D$ (cf. \cite{dgh,m4}). By the assumption, for
any codeword  $d\in D$, its components corresponding to $\om^+(\al)$
and $\om^-(\al)$ are the same. Therefore, the binary word
$(1100\cdots00)\in D^\perp$. Since $V^\natural$ is holomorphic,
$D^\perp=C$. Hence $C$ contains a codeword of weight $2$, which
contradicts $V^\natural_1=0$. Hence
$g=-1$ on $V_\Lambda^{T,+}$, which
completes the lemma.
\end{proof}

\begin{proposition}\label{prop:5.13} None of the Ising vectors in $V^\natural$ is of $\sigma$-type,
namely $V_{e}(\frac{1}{16})\neq0$ for any Ising vector $e\in V^\natural$.
\end{proposition}
\begin{proof}
Let $\alpha\in \Lambda(2)$. Then
$\tau_{\omega^+(\alpha)}\tau_{\omega^-(\alpha)}=z$ by the above
lemma. Now suppose that $e$ is of $\sigma$-type. By Lemma
\ref{LLSY}, $e\in
(V^{\natural})^{\langle\tau_{\omega^+(\alpha)},\tau_{\omega^-(\alpha)}\rangle}\subset
(V^{\natural})^z=V_\Lambda^+$ and hence $e$ is of $\sigma$-type in
$V_\Lambda^+$, which contradicts Corollary \ref{cor}.
\end{proof}

Proposition \ref{prp:5.11} and Lemma \ref{lem:5.12} show that the map $e\mapsto\tau_e$ from Ising vectors in $V_\Lambda^+\subset V^\natural$ to $\aut V^\natural$ is injective.

\begin{proposition}\label{prp:1to1} Let $e$, $f$ be Ising vectors in $V_\Lambda^+$ such that $\tau_e=\tau_f$ on $V^\natural$.
Then $e=f$.
\end{proposition}

\subsection{The centralizer of $\tau_e$ in $C_{\aut V^\natural}(z)$}

In this subsection, we shall determine the centralizer of $\tau_e$
in $C_{\aut V^\natural}(z)$ for any Ising vector $e$ in
$V_\Lambda^+$. Note that there are exactly $2$ types of Ising
vectors in $V_\Lambda^+$, namely, $\omega^\pm(\alpha)$ for some
$\alpha\in \Lambda(2)$ or $\varphi_x( e_E)$ for some sublattice
$E\cong \sqrt{2}E_8$ and $x\in \frac{1}{2}E$.

By a group of type $2^{n+m}$, we mean a group $H$ for which there
exists an exact sequence
\[
1 \to \Z_2^n \to H\to \Z_2^m \to 1
\]
which does not necessarily split. A group $G$ of type
$2^{n_1+\cdots+n_k}$ can be defined inductively by the exact
sequence
\[
1\to 2^{n_1+\cdots+n_{k-1}}\to G\to 2^{n_k}\to 1
\]

First, we consider Ising vectors in $V_\Lambda^+$ associated with
norm $4$ vectors in $\Lambda$.

\begin{lemma}\label{lem:gb1}
The centralizer of $z$ in $\aut V^\natural$ has the structure $2^{1+24}.Co_1$.
\end{lemma}
\begin{proof} Note that $V^\natural=V_\Lambda^+\oplus V_\Lambda^{T,+}$ is a simple current
extension\footnote{ An irreducible module $M$ of a simple VOA $U$ is
said to be a simple current module if the fusion product $M\times_U
W$ is also irreducible for all irreducible $U$-module $W$. A VOA $V$
is a simple current extension of a subVOA $U$ if $V$ is a direct sum
of simple current modules of $U$.}
 of $V_\Lambda^+$ and that for any automorphism $g\in \aut V_\Lambda^+$, the
$g$-conjugate module $g\circ V_\Lambda^{T,+}$, i.e, $g\circ
V_\Lambda^{T,+}=V_\Lambda^{T,+}$ as vector spaces and the vertex
operator $Y_g$ is defined by $Y_g(u,z)=Y(g^{-1}u,z)$, is isomorphic to
$V_\Lambda^{T,+}$ itself. In other word, $V_\Lambda^{T,+}$ is
preserved by $\aut V_\Lambda^+$.
Hence the following sequence is exact
(cf. \cite{sh}):
\begin{eqnarray*}
1\to\langle z\rangle\to C_{\aut V^\natural}(z)\overset{\mu}{\to}\aut
V_\Lambda^+\to1.
\end{eqnarray*}
\end{proof}

\begin{theorem} Let $\alpha\in\Lambda(2)$ and $\varepsilon\in\{\pm\}$.
Set $e=\omega^\varepsilon(\alpha)$. Then there is an exact sequence
\begin{eqnarray*}
1\to \{\varphi_\beta\in\Hom(\Lambda,\Z_2)\mid \langle
\beta,\alpha\rangle\in2\Z\}\to C_{\aut
V^\natural}(z,\tau_e)/\langle z\rangle\to Co_2\to1.
\end{eqnarray*}
In particular, the centralizer $C_{\aut V^\natural}(\tau_e,z)$ has
the structure $2^{2+22}.Co_2$.
\end{theorem}
\begin{proof} By Proposition \ref{lem:gb1}, $C_{\aut V^\natural}(z)/\langle z\rangle\cong\aut V_\Lambda^+$.
It follows from Proposition \ref{prp:1to1} that $C_{\aut
V^\natural}(z,\tau_e)/\langle z\rangle$ is isomorphic to the
stabilizer ${\rm Stab}_{\aut V_\Lambda^+}(e)$ of $e$ in $\aut
V_\Lambda^+$.

We view $\tau_e$ as an automorphism of $V_\Lambda^+$ and consider
the exact sequence in Theorem \ref{thm:4.10}:
\[
1 \longrightarrow  \mathrm{Hom}(\Lambda, \Z_2) \longrightarrow
\aut{V_\Lambda^+}\ \overset{\xi} \longrightarrow \aut
\Lambda/\langle\pm 1\rangle\longrightarrow  1.
\]
It follows from Proposition \ref{prp:3.2} that
$\tau_e=\varphi_\alpha\in\Hom(\Lambda,\Z_2)$ for some
$\alpha\in\Lambda$.

For $\varphi_\beta\in\Hom(\Lambda,\Z_2)$, $\beta\in\Lambda$, it is
easy to see that
$\varphi_\beta(\omega^\pm(\alpha))=\omega^\pm(\alpha)$ if $\langle
\alpha,\beta\rangle\in2\Z$, and
$\varphi_\beta(\omega^\pm(\alpha))=\omega^\mp(\alpha)$ if $\langle
\alpha,\beta\rangle\in2\Z+1$. Hence $\Hom(\Lambda,\Z_2)\cap {\rm
Stab}_{\aut V_\Lambda^+}(e)=\{\varphi_\beta\in\Hom(\Lambda,\Z_2)\mid
\langle \alpha,\beta\rangle\in2\Z\}\cong 2^{23}$.

Let $g$ be an element in $\aut V_\Lambda^+$ fixing $e$. Then
$g(\varphi_\alpha)=\varphi_{\xi(g)(\alpha)}$, hence $\xi(g)\in{\rm
Stab}_{Co_1}(\alpha+2\Lambda)=Co_2$. Conversely, let $p\in \aut
V_\Lambda^+$ be a lift of an element in $Co_2$. Then
$p(e)=\omega^+(\alpha)$ or $\omega^-(\alpha)$. If $p(e)\neq e$ then
$p\circ \varphi_\gamma(e)=e$, where $\gamma\in\Lambda$ such that
$\langle\gamma,\alpha\rangle\in2\Z+1$. This shows that $\xi({\rm
Stab}_{\aut V_\Lambda^+}(e))= Co_2$.

Thus we obtain an exact sequence
\begin{eqnarray*}
1\to \{\varphi_\beta\in\Hom(\Lambda,\Z_2)\mid \langle
\beta,\alpha\rangle\in2\Z\}\to {\rm Stab}_{\aut
V_\Lambda^+}(e)\to Co_2\to1.
\end{eqnarray*}
Since $\{\varphi_\beta\in\Hom(\Lambda,\Z_2)\mid \langle
\beta,\alpha\rangle\in2\Z\}=\langle\tau_{\omega^+(\beta)}\mid
\langle \beta,\alpha\rangle\in2\Z\rangle$ in $\aut V_\Lambda^+$, we
have
\[
\begin{split}
\mu^{-1}(\{\varphi_\beta\in\Hom(\Lambda,\Z_2)\mid \langle
\beta,\alpha\rangle\in2\Z\})&=\langle z, \tau_{\omega^+(\beta)}\mid
\langle \beta,\alpha\rangle\in2\Z\rangle\\
&\ \subset O_2(C_{\aut V^\natural}(z))\cong 2^{1+24}_+
\end{split}
\]
in $\aut V^\natural$, where $\mu$ is the natural epimorphism
$C_{\aut V^\natural}(z)\to C_{\aut V^\natural}(z)/\langle
z\rangle\cong\aut V_\Lambda^+$.

It is clear that $z$ and $\tau_e$ are in the center of $C_{\aut
V^\natural}(\tau_e,z)$. Hence
$\mu^{-1}(\{\varphi_\beta\in\Hom(\Lambda,\Z_2)\mid \langle
\beta,\alpha\rangle\in2\Z\})$ is of shape $2^{2+22}$, which proves
this theorem.
\end{proof}

If $e=\varphi_x(e_E)$, it comes out that the centralizer of $\tau_e$
in $C_{\aut V^\natural}(z)$ also stabilizes the VOA $V_E^+\subset
V_\Lambda^+$.

\begin{proposition}\label{prp:gb1} Let $E$ be a sublattice of $\Lambda$ isomorphic to $\sqrt2E_8$ and let $x\in \frac12E$.
Set $e=\varphi_x(e_E)$.
Then $C_{\aut
V_\Lambda^+}(\tau_{e})$ stabilizes $V_E^+$. Moreover,
$C_{\aut V^\natural}(\tau_{e},z)$ stabilizes $V_E^+$.
\end{proposition}
\begin{proof} Let $\mathcal{F}$ and $\mathcal{O}$ be an $8$-frame and an octad such that $E=E_{\mathcal{F}}(\mathcal{O})$.
Then by Proposition \ref{prp:4.17}, we obtain the exact sequence
$$1\to \{f_\alpha\in\Hom(\Lambda,\Z_2)\mid \langle\alpha,E\rangle\in2\Z\}\to C_{\aut V_\Lambda^+}(\tau_{e})\to C_{Co_1}(\varepsilon_{\mathcal{O}}^\mathcal{F})\to1. $$
Since $C_{Co_1}(\varepsilon_{\mathcal{O}}^\mathcal{F})$ preserves
$E$, $C_{\aut V_\Lambda^+}(\tau_{e})$ stabilizes $V_E^+$.
It follows from Lemma \ref{lem:gb1} and $z={\rm id}$ on
$V_\Lambda^+$ that $C_{\aut V^\natural}(\tau_{e},z)$
preserves $V_E^+$.
\end{proof}

In order to determine the centralizer $C_{\aut
V^\natural}(\tau_e,z)$ of $\tau_e$ and $z$ in $\aut V^\natural$, we
need the following lemma.

\begin{lemma}\label{prp:E} Let $E$ be a lattice isomorphic to $\sqrt2E_8$.
Then the stabilizer of $e_E$ in the subgroup $\aut
\hat{E}/\langle\theta\rangle$ of $\aut V_E^+$ is isomorphic to $\aut
E/\langle-1\rangle\cong O^+(8,2)$.
\end{lemma}
\begin{proof} Since $\langle E,E\rangle\subset 2\Z$, $\aut \hat{E}/\langle\theta\rangle$ is a split extension of $\aut E /\langle-1\rangle$ by $\Hom(E,\Z_2)$ (cf. \cite{flm}).
Clearly for $\varphi_\alpha\in\Hom(E,\Z_2)$, $\varphi_\alpha(e_E)=
e_E$ if and only if $\varphi_\alpha=1$. By the definition of $e_E$,
the subgroup $\aut E/\langle-1\rangle$ fixes $e_E$.
\end{proof}


\begin{theorem} Let $E$ be a sublattice of $\Lambda$ isomorphic to $\sqrt2E_8$ and $x$ a vector in $\frac{1}{2}E$.
Set $e=\varphi_x(e_E)$. Then the centralizer $C_{\aut
V^\natural}(\tau_e,z)$ has the structure $2^{2+8+16}.\Omega^+(8,2)$.
\end{theorem}
\begin{proof}
By Proposition \ref{prp:gb1}, $C_{\aut V^\natural}(\tau_e,z)/\langle
z\rangle$ is a subgroup of the stabilizer ${\rm Stab}_{\aut
V_\Lambda^+}(V_E^+)$ of $V_E^+$ in $\aut V_\Lambda^+$.

Let us describe ${\rm Stab}_{\aut V_\Lambda^+}(V_E^+)$. It follows
from Proposition \ref{prp:3.2} that $\Hom(\Lambda,\Z_2)$ is a
subgroup of ${\rm Stab}_{\aut V_\Lambda^+}(V_E^+)$. Let $g\in {\rm
Stab}_{\aut V_\Lambda^+}(V_E^+)$. Then $\xi(g)$ preserves $E$ in
$\Lambda$. Conversely, any element in $Co_1$ preserving $E$ lifts to
an element of ${\rm Stab}_{\aut V_\Lambda^+}(V_E^+)$. Hence
$$1\to\Hom(\Lambda,\Z_2)\to {\rm Stab}_{\aut V_\Lambda^+}(V_E^+)\to {\rm Stab}_{Co_1}(E)\to 1$$ is exact.
We note that ${\rm Stab}_{Co_1}(E)$ has the shape
$2^{1+8}\cdot\Omega^+(8,2)$ and that its subgroup acting on $E$ by
$\pm1$ is the maximal normal $2$-subgroup $O_2({\rm
Stab}_{Co_1}(E))\cong 2^{1+8}$. Clearly, $\tau_e$ is a lift of the
central element of ${\rm Stab}_{Co_1}(E)$ of order $2$.

Let $F$ be the subgroup of ${\rm Stab}_{\aut V_\Lambda^+}(V_E^+)$
acting trivially on $V_E^+$. Then the following sequence is exact:
$$1\to \{f_\alpha\in\Hom(\Lambda,\Z_2)\mid \langle \alpha,E\rangle\in2\Z\}\to F\to O_2({\rm Stab}_{Co_1}(E))\to 1.$$
In particular, $F\cong 2^{16}.2^{1+8}$. Consider the orthogonal
compliment $B$ of $E$ in $\Lambda$. Then $F$ acts on $V_B^+$ and,
the subgroup $Q$ of $F$ acting trivially on $V_B^+$ is described by
the following exact sequence
$$1\to \{f_\alpha\in\Hom(\Lambda,\Z_2)\mid \langle \alpha,B\oplus E\rangle\in2\Z\}\to Q\to Z({\rm Stab}_{Co_1}(E))\to1.$$
It is easy to see that $\{f_\alpha\in\Hom(\Lambda,\Z_2)\mid \langle
\alpha,B\oplus E\rangle\in2\Z\}=\{f_\alpha\mid \alpha\in
E\}\cong2^8$ and $Z({\rm Stab}_{Co_1}(E))\cong\Z_2$. Clearly,
$\tau_e\in Q$ and
$\tau_e\notin\Hom(\Lambda,\Z_2)$. Hence $Q\cong 2^9$.
Note that $\mu^{-1}(Q)=\langle z,\tau_{\omega^+(\alpha)},\tau_e\mid \alpha\in E\rangle\cong 2^{10}$.
Since $Q$ is normal in $F$, we obtain $F\cong 2^{9+16}$ and $\mu^{-1}(F)\cong 2^{10+16}$.

The group ${\rm Stab}_{\aut V_\Lambda^+}(V_E^+)/F\cong
2^8.\Omega^+(8,2)$ acts faithfully on $V_E^+$, and we view this
group as a subgroup of $\aut(\hat{E})/\langle\theta\rangle\cong 2^8:O^+(8,2)$. By
Lemma \ref{prp:E}, the stabilizer of $\tau_e$ in this group is
isomorphic to $\aut(E)^\prime/\langle-1\rangle\cong
\Omega^+(8,2)$. Thus we obtain the following exact sequence:
$$1\to F\to {\rm Stab}_{\aut V_\Lambda^+}(e)\to \aut(E)^\prime/\langle-1\rangle\to 1.$$
Hence $C_{\aut V^\natural}(\tau_e,z)\cong \mu^{-1}({\rm Stab}_{\aut V_\Lambda^+}(e))$ has the structure
$2^{10+16}.\Omega^+(8,2)$.
Since $z$ and $\tau_e$ are in the center, $C_{\aut V^\natural}(\tau_e,z)$ has the shape $2^{2+8+16}.\Omega^+(8,2)$
\end{proof}

\begin{remark} On the setting of the theorem above, $C_{\aut V^\natural}(\tau_e,z)$ acts on $V_E^+$ as $\Omega^+(8,2)$,
which is the quotient of the commutator subgroup of the Weyl group
of $E_8$ by its center. This result explains the $1A$ case of an
observation by Glauberman and Norton, which describes some
mysterious relations between the centralizer of $z$ and some $2A$
elements commuting $z$ in the Monster and the Weyl groups of certain
sublattices of the root lattice of type $E_8$ (cf. Section 4 of
\cite{gn}).
\end{remark}

\begin{remark}By \cite{co,m1}, an involution of the Monster is of type $2A$ if and only if it is a $\tau$-involution.
Hence if $e=\omega^\pm(\alpha)$, then $z$ and $\tau_e$ generate a
Klein's four group of type $(1A, 2A, 2A,2B)$ in the Monster because
$\langle
\tau_e,z\rangle=\{1,\tau_{\omega^+(\alpha)},\tau_{\omega^-(\alpha)},z\}$
by Lemma \ref{lem:5.12}. On the other hand, if $e=\varphi_x(e_E)$,
then $z$ and $\tau_e$ will generate a Klein's four group of type
$(1A, 2A, 2B, 2B)$ because $\langle \tau_e,z\rangle=\{1,\tau_e,
z\tau_e,z\}$ and $z\tau_e$ is not a $\tau$-involution by Proposition
\ref{prp:5.11}.
\end{remark}

\bibliographystyle{amsplain}

\end{document}